\documentclass[11pt,a4paper]{article}
\usepackage[T1]{fontenc}
\usepackage[utf8]{inputenc}
\usepackage{amsmath,amssymb,amsthm}
\usepackage{lmodern}
\usepackage[a4paper,left=26mm,right=26mm,top=25mm,bottom=25mm]{geometry}
\usepackage{microtype}
\usepackage{enumitem}
\usepackage{needspace}
\usepackage{tikz}
\usetikzlibrary{arrows.meta}
\usepackage[colorlinks=true,anchorcolor=blue,filecolor=blue,linkcolor=blue,urlcolor=blue,citecolor=blue]{hyperref}
\newtheorem{theorem}{Theorem}[section]
\newtheorem{lemma}[theorem]{Lemma}

\newtheorem{corollary}[theorem]{Corollary}
\newtheorem{conjecture}[theorem]{Conjecture}

\newtheorem{example}[theorem]{Example}
\newtheorem{remark}[theorem]{Remark}
\allowdisplaybreaks[1]
\newcommand{\mb}[1]{\text{\mathversion{bold}${#1}$}}

\begin{document}

\title{\Large\bfseries Proofs of
Li, Wang and Guo's conjectures on (signless) Laplacian spectral
radius\footnote{Partially supported by NSF of Fujian (No.\ 2026J001968).}}

\author{Hongzhang Chen$^{a}$, \, Jianxi Li$^{b}$\thanks{Corresponding
author.
Email addresses: mnhzchern@gmail.com,  ptjxli@hotmail.com, lxy09062021@126.com.}, \,  Xiaoyun Lv$^{b}$\\
$^a$School of Mathematics and Statistics,\\
Lanzhou University, Lanzhou, Gansu, P.R. China\\
$^b$School of Mathematics and Statistics,\\
Fujian Key Laboratory of
Granular Computing and Applications,\\
Minnan Normal University, Zhangzhou, Fujian, P.R. China.}

\date{}
\maketitle

\begin{abstract}
Let $\mu(G)$ and $q(G)$ be the Laplacian and signless Laplacian
spectral radii of a graph $G$. Li, Wang and Guo (2026) conjectured
that, if $H$ with $|V(H)|\geq 1$ is an induced subgraph of a graph
$G$, then $\mu(G)\le\max\limits_{v\in V(H)}\mu(H_v^t)$ and
$q(G)\le\max\limits_{v\in V(H)}q(H_v^t)$, where $t=e(G)-e(H)$ and
$H_v^t$ is the graph obtained from $H$ by adding $t$ pendant edges
to $v\in V(H)$. We prove this conjecture by establishing a spectral
comparison theorem for a broader family of positive semidefinite
graph matrices which includes the Laplacian and signless Laplacian
matrices. We also characterize all equality cases for the prescribed
subgraph $H$ with $e(H)\geq 1$ and $t\geq 1$.


\end{abstract}

\section{Introduction}\label{sec-1}

All graphs considered in this paper are finite, simple and undirected.
Let $G$ be a graph with vertex set $V(G)$ and edge set $E(G)$.
Its order and size are $|V(G)|$ and $e(G)=|E(G)|$, respectively.
For $v\in V(G)$, let $d_G(v)$ denote its degree, and let
$\Delta(G)$ denote the maximum degree of $G$.
A vertex is isolated if its degree is zero and pendant if
its degree is one.
For $U\subseteq V(G)$, we write $G[U]$ for the subgraph
induced by $U$. The complete graph of order $s$ is denoted by $K_s$.
We write $K_{r,s}$ for the complete bipartite graph with parts
of sizes $r$ and $s$.
In particular, $K_{1,s}$ is the star with $s$ edges.

Recall that the Laplacian and the signless Laplacian matrices of a
graph $G$ are $L(G)=D(G)-A(G)$ and $Q(G)=D(G)+A(G)$, where  $A(G)$
and $D(G)$ are the adjacency and the diagonal degree matrices of
$G$, respectively. We call the largest eigenvalues of $L(G)$ and
$Q(G)$ the Laplacian and the signless Laplacian spectral radii of
$G$, respectively, denoted by $\mu(G)$ and $q(G)$. Note that $L(G)$
and $Q(G)$ are positive semidefinite and have the same spectrum when
$G$ is bipartite. For every graph $G$, we have $\mu(G)\le q(G)$, and
the inequality is strict when $G$ is connected and nonbipartite. We
refer to Merris~\cite{Merris1994}, Brouwer and
Haemers~\cite{BH2012}, and Cvetkovi\'c, Rowlinson and
Simi\'c~\cite{CRS2007} for background on these matrices. Classical
bounds~\cite{AM1985,CRSBounds2007} involving vertex degrees and the
sums $d_G(u)+d_G(v)$ for $uv\in E(G)$ relate these spectral radii to
the local structure of $G$.

\begin{figure}[htb]
\centering
\begin{tikzpicture}[
  x=1.2cm,y=1cm,font=\small,
  vertex/.style={circle,fill=black,inner sep=0pt,minimum size=4pt},
  edge/.style={line width=0.65pt},
  extra/.style={edge,densely dashed},
  line cap=round,line join=round]
  \begin{scope}
    \filldraw[fill=black!5,draw=black!30,line width=0.4pt]
      (-0.12,0) ellipse [x radius=1.02,y radius=0.98];
    \coordinate (a) at (-0.72,0.52);
    \coordinate (b) at (-0.72,-0.52);
    \coordinate (v) at (0.56,0);
    \coordinate (x) at (2.05,0.85);
    \coordinate (y) at (2.05,-0.85);
    \draw[edge] (a)--(b)--(v)--(a);
    \draw[extra] (a)--(x)--(v)--(y)--(x);
    \foreach \u in {a,b,v,x,y} \node[vertex] at (\u) {};
    \node[below=3pt] at (v) {$v$};
    \node[inner sep=0pt] at (-0.12,0.79) {$H$};
    \node at (0.45,-1.5) {$G$};
  \end{scope}
  \draw[-{Stealth[length=2mm,width=1.4mm]},line width=0.65pt]
    (2.95,0)--(4.95,0)
    node[midway,above=5pt] {retain $H$}
    node[midway,below=5pt,font=\footnotesize] {attach $t$ new leaves to $v$};
  \begin{scope}[xshift=8.04cm]
    \filldraw[fill=black!5,draw=black!30,line width=0.4pt]
      (-0.12,0) ellipse [x radius=1.02,y radius=0.98];
    \coordinate (a) at (-0.72,0.52);
    \coordinate (b) at (-0.72,-0.52);
    \coordinate (v) at (0.56,0);
    \draw[edge] (a)--(b)--(v)--(a);
    \foreach \i/\px/\py in {1/2.1/1.1,2/2.4/0.4,3/2.4/-0.4,4/2.1/-1.1} {
      \coordinate (w\i) at (\px,\py);
      \draw[extra] (v)--(w\i);
      \node[vertex] at (w\i) {};
    }
    \foreach \u in {a,b,v} \node[vertex] at (\u) {};
    \node[below=3pt] at (v) {$v$};
    \node[inner sep=0pt] at (-0.12,0.79) {$H$};
    \node at (0.65,-1.5) {$H_v^t$};
  \end{scope}
\end{tikzpicture}
\caption{From $G$ to $H_v^t$ with $t=4$.}
\label{fig-hvt}
\end{figure}
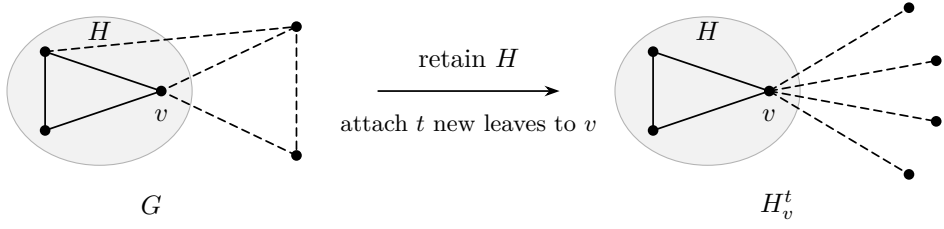

For $v\in V(H)$ and an integer $t\ge0$, let $H_v^t$
denote the graph obtained from $H$ by adding $t$ new vertices,
each adjacent only to $v$.
In particular, $H_v^0=H$.
This construction preserves $H$ as an induced subgraph
and makes all $t$ additional edges incident with $v$.
Figure~\ref{fig-hvt} illustrates the construction when
$t=e(G)-e(H)=4$.

With this notation, Li, Wang and Guo~\cite{LWG2026} proposed the following two conjectures.

\begin{conjecture}[\cite{LWG2026}, Conjecture 1]\label{conj-1.1}
Let $G$ contain $H$ as an induced subgraph, where $V(H)\ne\varnothing$,
and let $t=e(G)-e(H)$. Then $q(G)\le q(H_v^t)$ for some $v\in V(H)$.
\end{conjecture}

\begin{conjecture}[\cite{LWG2026}, Conjecture 2]\label{conj-1.2}
Let $G$ contain $H$ as an induced subgraph, where $V(H)\ne\varnothing$, and let $t=e(G)-e(H)$. Then $\mu(G)\le\mu(H_v^t)$ for some $v\in V(H)$.
\end{conjecture}

Li, Wang and Guo~\cite[Theorem~18]{LWG2026} proved
Conjecture~\ref{conj-1.1} under the condition
$e(G)\ge3(e(H)-\Delta(H))+3$.
They also proved both conjectures when $H$ is a cycle,
a path or a clique~\cite[Section~4]{LWG2026}.
Both conjectures allow arbitrary $H$ and every number
of additional edges.
In particular, a component of $G-V(H)$ may have neighbors
at several vertices of $H$ and may itself contain cycles.


 In this paper, we prove Conjectures~\ref{conj-1.1}
and~\ref{conj-1.2} for every induced subgraph $H$ with $|V(H)|\geq
1$ and every $t\ge0$ by establishing a spectral comparison theorem
for a broader family of graph matrices $M_\alpha(G)=D(G)+\alpha
A(G)$, where $-1\le\alpha\le1$. Let $\rho_\alpha(G)$ be the spectral
radius of $M_\alpha(G)$.  Our main result is the following.


\begin{theorem}\label{thm-1.3}
Let $H$ with $|V(H)|\geq 1$ be an induced subgraph of a graph $G$.
Then for $\alpha\in [-1, 1]$, we have
\begin{equation}\label{e-1}
\rho_\alpha(G)\le\max_{v\in V(H)}\rho_\alpha(H_v^t),
\end{equation}
 where $t=e(G)-e(H)$.
\end{theorem}

Since $\rho_{-1}(G)=\mu(G)$ and $\rho_1(G)=q(G)$, we obtain
the following corollary by taking $\alpha=-1$ and $\alpha=1$
in Theorem~\ref{thm-1.3}. This proves
Conjectures~\ref{conj-1.1} and~\ref{conj-1.2}.

\begin{corollary}\label{cor-1.4}
Let $H$ with $|V(H)|\geq 1$ be an induced subgraph of a graph $G$.
Then we have
\begin{gather*}
\mu(G)\le\max_{v\in V(H)}\mu(H_v^t)\, \mbox{ and }\,
q(G)\le\max_{v\in V(H)}q(H_v^t),
\end{gather*}
where $t=e(G)-e(H)$.
\end{corollary}

Note that each graph $H_v^t$ has $e(H)+t$ edges and contains $H$ as an induced
subgraph. 
Thus~\eqref{e-1} determines the exact maximum of
$\rho_\alpha$ among all graphs satisfying these two conditions. In fact, a maximizing attachment vertex may depend on
$\alpha$ and $t$. The following theorem characterizes the equality cases relative to
the prescribed induced copy of $H$. In each case, any
number of isolated vertices may be added outside $H$.

\begin{theorem}\label{thm-1.5}
For $|\alpha|\le1$, let $H$ with $e(H)\ge1$ be an induced
subgraph of a graph $G$. Suppose that $t=e(G)-e(H)\ge1$. Then the equality holds in~\eqref{e-1} if and only if one of the following
conditions holds.
\begin{enumerate}[label=\textup{(\roman*)},leftmargin=2em]
\item
The graph $G$ is obtained from $H_v^t$ by adding isolated
vertices outside $H$, where $v\in V(H)$ satisfies
$\rho_\alpha(H_v^t)=\max_{u\in V(H)}\rho_\alpha(H_u^t)$.
\item
We have $\alpha=1$, $t=2$ and $\Delta(H)=1$, and $G$
is obtained from $H$ by adding one vertex adjacent to both
endpoints of a $K_2$ component, together with isolated
vertices outside $H$.
\end{enumerate}
\end{theorem}

\begin{remark}\label{rem-1.6}{\rm
The condition $e(H)\ge1$ in Theorem~\ref{thm-1.5} cannot be omitted. For example, if $e(H)=0$ and $t\ge1$, then each graph $H_v^t$
consists of $K_{1,t}$ and isolated vertices, so its spectral radius
is given by~\eqref{e-4} with $s=t$. Moreover, the equality in~\eqref{e-1} may hold even when all additional edges are disjoint from the prescribed
copy of $H$. For example, let $G=K_2\cup K_1$, and take $H$ to be
the subgraph induced by its isolated vertex $v$. Then $t=1$, and
$\rho_\alpha(G)=\rho_\alpha(H_v^1)=1+|\alpha|$ for every
$|\alpha|\le1$, although neither endpoint of the unique edge of $G$
belongs to $H$.}
\end{remark}

Section~\ref{sec-3} proves the comparison theorem and characterizes
equality. Section~\ref{sec-5} presents a counterexample to the adjacency
analogue and establishes the comparison under an additional
condition on the components of $G-V(H)$.

\section{Preliminaries}\label{sec-2}

For a real symmetric matrix $X$, let $\lambda_{\max}(X)$ be the
 largest eigenvalue of $X$. We write $X\succeq0$ if $X$ is positive
semidefinite. For real symmetric matrices $X$ and $Y$ of the same
order, $X\preceq Y$ means $Y-X\succeq0$. We write $I$ for the identity matrix. Let $\mb e_v$ be the standard coordinate vector with
entry $1$ at $v$ and zeros elsewhere. We use $\oplus$ for the direct
sum of matrices and $\|\mb x\|=(\mb x^{\mathsf T}\mb x)^{1/2}$ for
the Euclidean norm of a real vector $\mb x$. Matrix and vector
dimensions are determined by the context unless otherwise specified.

We use the following standard matrix results from
Horn and Johnson~\cite{HJ2013}.

\begin{lemma}[\cite{HJ2013}]\label{lem-2.1}
Let $W$ be a real matrix. The nonzero eigenvalues of $W^{\mathsf T}W$ and $WW^{\mathsf T}$ coincide, including multiplicities. In particular, $W^{\mathsf T}W\preceq I$ if and only if $WW^{\mathsf T}\preceq I$, with identity matrices of the appropriate orders.
\end{lemma}

\begin{lemma}[\cite{HJ2013}]\label{lem-2.2}
Let $X$ be real symmetric and positive definite, let $Y$ be real symmetric, and let $Z$ have the appropriate dimensions. Then
\[
\begin{pmatrix}Y&Z\\Z^{\mathsf T}&X\end{pmatrix}\succeq0
\quad\Longleftrightarrow\quad
Y-ZX^{-1}Z^{\mathsf T}\succeq0.
\]
The block matrix is singular if and only if its Schur complement $Y-ZX^{-1}Z^{\mathsf T}$ is singular.
\end{lemma}

We also use the variational characterization of the largest
eigenvalue~\cite{HJ2013}. If $X\succeq0$ and $\mb x^{\mathsf T}X \mb
x=0$, then $X \mb x=0$. This follows by expanding $\mb x$ in an
orthonormal eigenbasis of $X$.

For a graph $F$ and a vector $\mb x\in\mathbb{R}^{|V(F)|}$, we have
\begin{equation}\label{e-2}
\begin{aligned}
\mb x^{\mathsf T}M_\alpha(F) \mb x
&=\sum_{uv\in E(F)}(x_u^2+x_v^2+2\alpha x_ux_v)\\
&=\sum_{uv\in E(F)} \bigl[(x_u+\alpha
x_v)^2+(1-\alpha^2)x_v^2\bigr],
\end{aligned}
\end{equation}
where the endpoints of each edge may be ordered arbitrarily.
For $|\alpha|\le1$, each summand is nonnegative, so
$M_\alpha(F)\succeq0$.
By~\eqref{e-2}, for every subgraph $F$ of a graph $G$,
\[
M_\alpha(F)\oplus0\preceq M_\alpha(G),
\]
where the zero block is indexed by $V(G)\setminus V(F)$.
Hence $\rho_\alpha(F)\le\rho_\alpha(G)$ by the variational
characterization of the largest eigenvalue.

Moreover, $2\alpha x_ux_v\ge-|\alpha|(x_u^2+x_v^2)$ implies
\begin{equation}\label{e-3}
\begin{aligned}
\mb x^{\mathsf T}M_\alpha(F) \mb x
&\ge(1-|\alpha|)\sum_{uv\in E(F)}(x_u^2+x_v^2)\\
&=(1-|\alpha|)\sum_{v\in V(F)}d_F(v)x_v^2.
\end{aligned}
\end{equation}
For $|\alpha|<1$, by~\eqref{e-3} and the zero rows and columns
at isolated vertices,
\[
\ker M_\alpha(F)=\ker D(F)
=\{\mb x:x_v=0\text{ whenever }d_F(v)>0\}.
\]

For $\alpha=-1$ and $\alpha=1$, respectively,
identity~\eqref{e-2} reduces to $\mb x^{\mathsf T}L(F)\mb
x=\sum_{uv\in E(F)}(x_u-x_v)^2$ and $\mb x^{\mathsf T}Q(F)\mb
x=\sum_{uv\in E(F)}(x_u+x_v)^2$.
\begin{lemma}\label{lem-2.3}
For every integer $s\ge1$ and $|\alpha|\le1$, we have
\begin{equation}\label{e-4}
\rho_\alpha(K_{1,s})=
\frac{s+1+\sqrt{(s-1)^2+4\alpha^2s}}2.
\end{equation}
Moreover, $\rho_\alpha(K_{1,s}) > \max\{s,1\}$ for
$\alpha\ne0$.
\end{lemma}

\begin{proof}
Index the center of $K_{1,s}$ by $0$ and its leaves by $1,\ldots,s$.
If $x_0=0$ and $\sum_{i=1}^s x_i=0$, then $M_\alpha(K_{1,s}) \mb x=
\mb x$. These vectors form an invariant subspace of dimension $s-1$.
Its orthogonal complement is
$\operatorname{span}\{\mb e_0,\mb e_1+\cdots+\mb e_s\}$, which is also invariant
under $M_\alpha(K_{1,s})$. With respect to the ordered basis
$(\mb e_0,\mb e_1+\cdots+\mb e_s)$, the restriction of $M_\alpha(K_{1,s})$ to
this complement is represented by
\begin{equation*}
\begin{pmatrix}s&\alpha s\\\alpha&1\end{pmatrix}.
\end{equation*}
The characteristic polynomial of this matrix is
\begin{equation*}
(z-s)(z-1)-\alpha^2s
=z^2-(s+1)z+s(1-\alpha^2).
\end{equation*}
The larger root is the expression in~\eqref{e-4}.
Since $\sqrt{(s-1)^2+4\alpha^2s}\ge|s-1|$, it is at least
$\max\{s,1\}$, with strict inequality when $\alpha\ne0$.
The remaining $s-1$ eigenvalues of $M_\alpha(K_{1,s})$ are $1$,
so this root equals $\rho_\alpha(K_{1,s})$.
\end{proof}

\section{Proofs of Theorems~\ref{thm-1.3} and~\ref{thm-1.5}}\label{sec-3}


We first prove the spectral comparison in Theorem~\ref{thm-1.3}.
We then use the equality conditions in its proof to establish
Theorem~\ref{thm-1.5}.

\begin{proof}[Proof of Theorem~\ref{thm-1.3}]
If $t=0$, then $G$ is obtained from $H$ by adding isolated
vertices, so equality holds. In what follows, we assume that
$t\ge1$. We now consider the following two cases.

\vspace{1mm}
\noindent\textbf{Case 1.} $\alpha=0$.
\vspace{1mm}

We have $M_0(G)=D(G)$ and hence
$\rho_0(G)=\Delta(G)$. Since $|E(G)\setminus E(H)|=t$, we have
$d_G(u)\le d_H(u)+t\le\Delta(H)+t$ for $u\in V(H)$, whereas
$d_G(u)\le t$ for $u\in V(G)\setminus V(H)$. Moreover,
$\Delta(H_v^t)\le\Delta(H)+t$ for every $v\in V(H)$, with equality
whenever $d_H(v)=\Delta(H)$. Thus
\begin{equation*}
\rho_0(G)=\Delta(G)\le\Delta(H)+t =\max_{v\in V(H)}\rho_0(H_v^t),
\end{equation*}
as desired.

\vspace{1mm}
\noindent\textbf{Case 2.} $\alpha\neq 0$, that is $0<|\alpha|\le1$.
\vspace{1mm}

Let $\lambda=\max_{v\in V(H)}\rho_\alpha(H_v^t)$. It suffices to show that $\lambda I-M_\alpha(G)\succeq0$.


Choose a unit vector $\mb x$ satisfying
$M_\alpha(H)\mb x=\rho_\alpha(H)\mb x$, and take $v\in V(H)$
with $x_v\ne0$. By the variational characterization applied
to the extension of $\mb x$ by zero at the new leaves,
\[
\lambda\ge\rho_\alpha(H_v^t)
\ge\rho_\alpha(H)+t x_v^2>\rho_\alpha(H).
\]
Also, $K_{1,t}$ is a subgraph of every $H_v^t$, so
Lemma~\ref{lem-2.3} and subgraph monotonicity imply
\[
\lambda\ge\rho_\alpha(K_{1,t})>1.
\]

With the vertices of $H$ listed first, define
$R=[\lambda I-(M_\alpha(H)\oplus0)]^{-1}$,
where the zero block has order $|V(G)|-|V(H)|$.
The preceding inequalities show that $R$ is well defined and
positive definite.
Let $p=(\lambda-1)/[t(\lambda-1+\alpha^2)]-\lambda^{-1}$.

\vspace{3mm}
\noindent\textbf{Claim 1.}
\emph{The matrix $R-\lambda^{-1}I$ is positive semidefinite and
has kernel $\ker(M_\alpha(H)\oplus0)$.
Moreover, $p\ge0$ and $|(R-\lambda^{-1}I)_{uv}|\le p$
for all $u,v\in V(H)$.}

\begin{proof}[Proof of Claim 1.]
Each eigenvalue $\theta$ of $M_\alpha(H)\oplus0$ satisfies
$0\le\theta<\lambda$ and corresponds to the eigenvalue
$(\lambda-\theta)^{-1}-\lambda^{-1}
=\theta/[\lambda(\lambda-\theta)]$
of $R-\lambda^{-1}I$.
This value is nonnegative and vanishes if and only if $\theta=0$.
Hence
\begin{equation}\label{e-7}
\begin{gathered}
R-\frac1\lambda I\succeq0\, \mbox{ and } \,
\ker\left(R-\frac1\lambda I\right)
=\ker(M_\alpha(H)\oplus0).
\end{gathered}
\end{equation}
Also, $R=(\lambda I-M_\alpha(H))^{-1}\oplus\lambda^{-1}I$,
so all entries of $R-\lambda^{-1}I$ outside its principal block
indexed by $V(H)$ are zero.

Fix $v\in V(H)$.
Since $\lambda\ge\rho_\alpha(H_v^t)$, we have
\begin{equation*}
\lambda I-M_\alpha(H_v^t)=
\begin{pmatrix}
\lambda I-M_\alpha(H)-t \mb e_v \mb e_v^{\mathsf T}
&-\alpha \mb e_v \mathbf1^{\mathsf T}\\
-\alpha\mathbf1 \mb e_v^{\mathsf T} &(\lambda-1)I_t
\end{pmatrix}\succeq0,
\end{equation*}
where $\mathbf1\in\mathbb R^t$ is the vector of ones.
Since $\lambda>1$, Lemma~\ref{lem-2.2} shows that the
Schur complement of $(\lambda-1)I_t$ satisfies
\begin{equation*}
\lambda I-M_\alpha(H) -t\left(1+\frac{\alpha^2}{\lambda-1}\right)
\mb e_v \mb e_v^{\mathsf T}\succeq0.
\end{equation*}
Evaluating its quadratic form at
$(\lambda I-M_\alpha(H))^{-1}\mb e_v$, we obtain
\[
R_{vv}-t\left(1+\frac{\alpha^2}{\lambda-1}\right)R_{vv}^2\ge0.
\]
Since $R_{vv}\ge\lambda^{-1}>0$ by~\eqref{e-7}, we obtain
\begin{equation}\label{e-9}
\frac1\lambda\le R_{vv}\le
\frac{\lambda-1}{t(\lambda-1+\alpha^2)}.
\end{equation}
By~\eqref{e-9}, we have $0\le R_{vv}-\lambda^{-1}\le p$
for every $v\in V(H)$, so $p\ge0$.
The nonnegative $2\times2$ principal minors of
$R-\lambda^{-1}I$ imply, for all $u,v\in V(H)$,
\begin{equation}\label{e-11}
\left|\left(R-\frac1\lambda I\right)_{uv}\right|^2
\le
\left(R_{uu}-\frac1\lambda\right)
\left(R_{vv}-\frac1\lambda\right)
\le p^2.
\end{equation}
This proves Claim~1.
\end{proof}

\smallskip
Write $E(G)\setminus E(H)=\{e_1,\ldots,e_t\}$, where $e_i=u_iv_i$.
Since $H$ is induced, each $e_i$ has at most one endpoint in $V(H)$.
We may therefore label its endpoints so that $v_i\notin V(H)$ for
every $i$. For each $1\le i\le t$, let $B_i$ be the matrix with
columns $\mb e_{u_i}+\alpha \mb e_{v_i}$ and $\sqrt{1-\alpha^2}\,
\mb e_{v_i}$.
Write $B=(B_1\ \cdots\ B_t)$.
By~\eqref{e-2}, we have
\begin{equation}\label{e-12}
M_\alpha(G)=M_\alpha(H)\oplus0+BB^{\mathsf T}.
\end{equation}
For $\alpha=\pm1$, the second column of each $B_i$ is zero,
and the same factorization remains valid.

Let $\mb y=(\mb y_1^{\mathsf T},\ldots, \mb y_t^{\mathsf
T})^{\mathsf T} \in \mathbb R^{2t}$, with $\mb y_i\in\mathbb R^2$
for each $i$. Set $\mb f=(1,0)^{\mathsf T}$ and $\mb
g=(\alpha,\sqrt{1-\alpha^2})^{\mathsf T}$. Then $\|\mb f\|=\|\mb
g\|=1$, $\mb f^{\mathsf T}\mb g=\alpha$, and $B_i \mb y_i=(\mb
f^{\mathsf T}\mb y_i)\mb e_{u_i} +(\mb g^{\mathsf T}\mb y_i)\mb
e_{v_i}$. For each $i$, define $s_i=\max\{|\mb f^{\mathsf T}\mb
y_i|,\, |\mb g^{\mathsf T}\mb y_i|\}$ and $z_i=\min\{|\mb f^{\mathsf
T}\mb y_i|,\, |\mb g^{\mathsf T}\mb y_i|\}$. Thus $\|B_i \mb
y_i\|^2=s_i^2+z_i^2$.

\vspace{3mm}
\noindent\textbf{Claim 2.}
\emph{For every $\mb y\in\mathbb R^{2t}$, we have
$\mb y^{\mathsf T}B^{\mathsf T}RB\mb y
\le\sum_{i=1}^t[(\lambda-1)s_i^2/(\lambda-1+\alpha^2)+z_i^2/\lambda]$.}

\begin{proof}[Proof of Claim 2.]
For $i\ne j$, the edges $e_i$ and $e_j$ share at most one endpoint,
so $\langle B_i \mb y_i, B_j \mb y_j\rangle\le s_is_j$.
Consequently,
\begin{equation}\label{e-13}
\begin{aligned}
\|B \mb y\|^2 &=\sum_{i=1}^t\|B_i \mb y_i\|^2
  +2\sum_{1\le i<j\le t}\langle B_i \mb y_i, B_j \mb y_j\rangle\\
&\le\sum_{i=1}^t(s_i^2+z_i^2)
  +2\sum_{1\le i<j\le t}s_is_j\\
&=\left(\sum_{i=1}^t s_i\right)^2+\sum_{i=1}^t z_i^2.
\end{aligned}
\end{equation}
Since $v_i\notin V(H)$ for every $i$ and
$R-\lambda^{-1}I$ is zero outside its principal block
indexed by $V(H)$, we obtain from~\eqref{e-11}
\begin{equation}\label{e-14}
\begin{aligned}
\mb y^{\mathsf T}B^{\mathsf T}\left(R-\frac1\lambda I\right)B \mb y
&=\sum_{\substack{1\le i,j\le t\\u_i,u_j\in V(H)}} (\mb f^{\mathsf
T}\mb y_i)(\mb f^{\mathsf T}\mb y_j)
\left(R-\frac1\lambda I\right)_{u_i u_j}\\
&\le p\left( \sum_{\substack{1\le i\le t\\u_i\in V(H)}} |\mb
f^{\mathsf T} \mb y_i|
\right)^2
\le p\left(\sum_{i=1}^t s_i\right)^2.
\end{aligned}
\end{equation}
Combining~\eqref{e-13} and~\eqref{e-14} and applying
Cauchy's inequality, we obtain
\begin{equation}\label{e-15}
\begin{aligned}
\mb y^{\mathsf T}B^{\mathsf T}RB \mb y
&\le\left(\frac1\lambda+p\right) \left(\sum_{i=1}^t s_i\right)^2
+\frac1\lambda\sum_{i=1}^t z_i^2\\
&\le t\left(\frac1\lambda+p\right)\sum_{i=1}^t s_i^2
+\frac1\lambda\sum_{i=1}^t z_i^2\\
&=\sum_{i=1}^t\left(
\frac{\lambda-1}{\lambda-1+\alpha^2}s_i^2
+\frac1\lambda z_i^2\right),
\end{aligned}
\end{equation}
as desired.
\end{proof}

\vspace{3mm}
\noindent\textbf{Claim 3.}
\emph{For every $1\le i\le t$, the $i$th summand in the last
line of~\eqref{e-15} is at most $\|\mb y_i\|^2$.}

\begin{proof}[Proof of Claim 3.]
Consider the symmetric matrix of order two
\begin{equation*}
\frac{\lambda-1}{\lambda-1+\alpha^2}\mb f \mb f^{\mathsf T}
+\frac1\lambda \mb g \mb g^{\mathsf T}.
\end{equation*}
Since $\|\mb f\|=\|\mb g\|=1$ and
$\mb f^{\mathsf T}\mb g=\alpha$, its characteristic polynomial is
\[
(z-1)\left(
z-\frac{(\lambda-1)(1-\alpha^2)}
        {\lambda(\lambda-1+\alpha^2)}
\right).
\]
The second root lies in $[0,\lambda^{-1}]\subset[0,1)$,
so the largest eigenvalue is $1$.
The same calculation applies after interchanging $\mb f$ and
$\mb g$. By the definitions of $s_i$ and $z_i$, the $i$th
summand in~\eqref{e-15} is the quadratic form at $\mb y_i$
of one of these two matrices. Hence
\begin{equation}\label{e-17}
\frac{\lambda-1}{\lambda-1+\alpha^2}s_i^2 +\frac1\lambda z_i^2
\le\|\mb y_i\|^2.
\end{equation}
This proves Claim~3.
\end{proof}

\smallskip
By Claims~2 and~3, we have $B^{\mathsf T}RB\preceq I$.
Applying Lemma~\ref{lem-2.1} to $R^{1/2}B$, we obtain
$R^{1/2}BB^{\mathsf T}R^{1/2}\preceq I$.
Thus, by~\eqref{e-12} and the definition of $R$,
\begin{equation}\label{e-18}
\begin{aligned}
\lambda I-M_\alpha(G)
&=R^{-1}-BB^{\mathsf T}\\
&=R^{-1/2}
\left(I-R^{1/2}BB^{\mathsf T}R^{1/2}\right)
R^{-1/2}\\
&\succeq0.
\end{aligned}
\end{equation}
Thus $\rho_\alpha(G)\le\lambda$, completing the proof.
\end{proof}

\begin{proof}[Proof of Theorem~\ref{thm-1.5}]
\noindent\textbf{Necessity.}
Suppose that equality holds in~\eqref{e-1}.

\vspace{1mm}
\noindent\textbf{Case 1.} $\alpha=0$.
\vspace{1mm}

In this case, equality is equivalent to $\Delta(G)=\Delta(H)+t$.
Since $\Delta(H)\ge1$, every $u\in V(G)\setminus V(H)$ satisfies
$d_G(u)\le t<\Delta(H)+t$. Hence some $v\in V(H)$ satisfies
$d_G(v)=\Delta(H)+t$. The inequalities $d_G(v)\le
d_H(v)+t\le\Delta(H)+t$ then imply $d_H(v)=\Delta(H)$ and
$d_G(v)-d_H(v)=t$. Thus every additional edge is incident with $v$.
Since $H$ is induced, these edges join $v$ to $t$ distinct vertices
outside $H$, each of degree one in $G$. Every other vertex outside
$H$ is isolated.
This is the configuration in part~\textup{(i)}, with
$d_H(v)=\Delta(H)$.

\vspace{1mm}
\noindent\textbf{Case 2.} $\alpha\ne0$.
\vspace{1mm}

Retain the notation from the proof of Theorem~\ref{thm-1.3}.
Since $e(H)\ge1$,~\eqref{e-7} shows that
$R-\lambda^{-1}I$ is nonzero and positive semidefinite.
Its nonzero entries lie in the principal block indexed by $V(H)$,
so~\eqref{e-9} implies
\[
p\ge\max_{v\in V(H)}(R_{vv}-\lambda^{-1})>0.
\]
By~\eqref{e-18} and Lemma~\ref{lem-2.1},
$\rho_\alpha(G)=\lambda$ implies that $1$ is an eigenvalue of
$B^{\mathsf T}RB$.
Choose a unit vector $\mb y\in\mathbb R^{2t}$ satisfying
$B^{\mathsf T}RB\mb y=\mb y$, and define $s_i$ and $z_i$
as in the proof of Theorem~\ref{thm-1.3}.
Since $\mb y^{\mathsf T}B^{\mathsf T}RB\mb y=\|\mb y\|^2=1$,
equality holds throughout~\eqref{e-13}, \eqref{e-14},
\eqref{e-15} and~\eqref{e-17}.
Equality in Cauchy's inequality implies $s_1=\cdots=s_t=s$.
Moreover, $B\mb y\ne0$, so $s>0$.

\vspace{3mm}
\noindent\textbf{Claim 1.}
\emph{For every $1\le i\le t$, the additional edge $e_i$ joins
a nonisolated vertex $u_i\in V(H)$ to a vertex $v_i\notin V(H)$.
Moreover, $|\mb f^{\mathsf T}\mb y_i|=s$ and
$R_{u_i u_i}=\lambda^{-1}+p$.}

\begin{proof}[Proof of Claim 1.]
Since $p>0$, equality in the last inequality of~\eqref{e-14}
implies
\[
\sum_{u_i\in V(H)}|\mb f^{\mathsf T}\mb y_i|=ts.
\]
Each summand is at most $s>0$, so $u_i\in V(H)$ and
$|\mb f^{\mathsf T}\mb y_i|=s$ for every $i$.
Equality in the first inequality of~\eqref{e-14} forces each
diagonal term to attain its upper bound. Hence
$R_{u_i u_i}-\lambda^{-1}=p>0$.
Since $R_{uu}=\lambda^{-1}$ for every isolated vertex $u$ of $H$,
we obtain $d_H(u_i)>0$ for every $i$.
\end{proof}

\vspace{3mm}
\noindent\textbf{Claim 2.}
\emph{All additional edges have a common endpoint.}

\begin{proof}[Proof of Claim 2.]
If two additional edges $e_i$ and $e_j$ were disjoint, then $\langle
B_i \mb y_i, B_j \mb y_j\rangle=0<s^2$, contradicting equality
in~\eqref{e-13}. Thus every two additional edges intersect. These
edges form a bipartite graph with parts $V(H)$ and $V(G)\setminus
V(H)$. If $t\ge2$, choose two distinct additional edges $uv$ and
$uw$. Any additional edge not incident with $u$ would have to be
$vw$ to intersect both, which is impossible because $v$ and $w$
belong to the same part. Hence all additional edges have a common
endpoint and form a star. The same conclusion holds when $t=1$.
This proves Claim~2.
\end{proof}

By Claim~2, the additional edges form a star.
We distinguish two subcases according to the position of its
center. When $t=1$, we choose its endpoint in $H$ as the center.
Figure~\ref{fig-2} illustrates these positions and the exceptional
configuration in part~\textup{(ii)}.
The internal edges of $H$ are omitted in panels \textup{(a)} and
\textup{(b)}, and panel \textup{(c)} shows only the affected component.

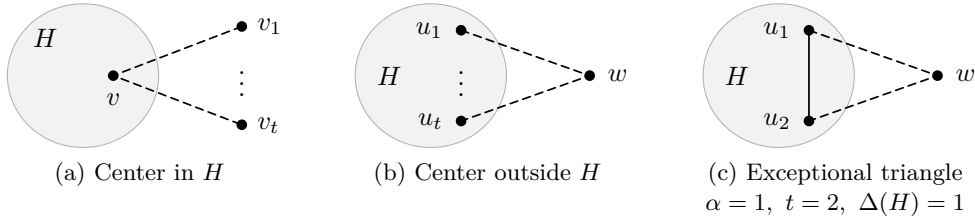
\begin{figure}[htb]
\centering
\begin{tikzpicture}[
  x=1cm,y=1cm,font=\small,
  vertex/.style={circle,fill=black,inner sep=0pt,minimum size=4pt},
  edge/.style={line width=0.65pt},
  extra/.style={edge,densely dashed},
  line cap=round,line join=round]
  \begin{scope}
    \filldraw[fill=black!5,draw=black!30,line width=0.4pt]
      (0,0) ellipse [x radius=1,y radius=0.95];
    \coordinate (v) at (0.4,0);
    \coordinate (vone) at (2.1,0.65);
    \coordinate (vt) at (2.1,-0.65);
    \draw[extra] (vone)--(v)--(vt);
    \foreach \u in {v,vone,vt} \node[vertex] at (\u) {};
    \node[below=3pt] at (v) {$v$};
    \node[right=2pt] at (vone) {$v_1$};
    \node[right=2pt] at (vt) {$v_t$};
    \node at (2.1,0) {$\vdots$};
    \node at (-0.5,0.5) {$H$};
    \node[font=\footnotesize] at (0.75,-1.3)
      {\textup{(a)} Center in $H$};
  \end{scope}
  \begin{scope}[xshift=4.6cm]
    \filldraw[fill=black!5,draw=black!30,line width=0.4pt]
      (0,0) ellipse [x radius=1,y radius=0.95];
    \coordinate (uone) at (0.4,0.6);
    \coordinate (ut) at (0.4,-0.6);
    \coordinate (w) at (2.1,0);
    \draw[extra] (uone)--(w)--(ut);
    \foreach \u in {uone,ut,w} \node[vertex] at (\u) {};
    \node[left=3pt] at (uone) {$u_1$};
    \node[left=3pt] at (ut) {$u_t$};
    \node[right=3pt] at (w) {$w$};
    \node at (0.4,0) {$\vdots$};
    \node at (-0.55,0) {$H$};
    \node[font=\footnotesize] at (0.75,-1.3)
      {\textup{(b)} Center outside $H$};
  \end{scope}
  \begin{scope}[xshift=9.2cm]
    \filldraw[fill=black!5,draw=black!30,line width=0.4pt]
      (0,0) ellipse [x radius=1,y radius=0.95];
    \coordinate (uone) at (0.4,0.6);
    \coordinate (utwo) at (0.4,-0.6);
    \coordinate (w) at (2.1,0);
    \draw[edge] (uone)--(utwo);
    \draw[extra] (uone)--(w)--(utwo);
    \foreach \u in {uone,utwo,w} \node[vertex] at (\u) {};
    \node[left=3pt] at (uone) {$u_1$};
    \node[left=3pt] at (utwo) {$u_2$};
    \node[right=3pt] at (w) {$w$};
    \node at (-0.55,0) {$H$};
    \node[font=\footnotesize] at (0.75,-1.3)
      {\textup{(c)} Exceptional triangle};
    \node[font=\footnotesize] at (0.75,-1.7)
      {$\alpha=1,\ t=2,\ \Delta(H)=1$};
  \end{scope}
\end{tikzpicture}
\caption{Star configurations. Dashed edges are added to $H$.}
\label{fig-2}
\end{figure}

\vspace{1mm}
\noindent\textbf{Subcase 2.1.} \emph{The center belongs to $H$.}
\vspace{1mm}

Let $v\in V(H)$ be the center. Every other endpoint of an
additional edge is pendant in $G$, and all remaining vertices
outside $H$ are isolated.
Thus $G$ is obtained from $H_v^t$ by adding isolated vertices.
Since $\rho_\alpha(H_v^t)=\rho_\alpha(G)=\lambda$,
part~\textup{(i)} follows.

\vspace{1mm}
\noindent\textbf{Subcase 2.2.} \emph{The center lies outside $H$.}
\vspace{1mm}

Here $t\ge2$. Let $u_1,\ldots,u_t$ be its distinct neighbors in $H$.
We show that this configuration forces part~\textup{(ii)}.

\vspace{3mm}
\noindent\textbf{Claim 3.}
\emph{$|\alpha|=1$. Moreover,
$M_\alpha(H)(\mb e_{u_i}-\mb e_{u_j})=0$
for every pair of distinct indices $i,j$.}

\begin{proof}[Proof of Claim 3.]
Write $\mb f^{\mathsf T}\mb y_i=\varepsilon_i s$, where
$\varepsilon_i\in\{-1,1\}$. In~\eqref{e-14}, each term of the
double sum is at most $ps^2$, while the sum is $pt^2s^2$.
Thus, for all $1\le i,j\le t$,
\[
\left(R-\frac1\lambda I\right)_{u_i u_j}
=p\varepsilon_i\varepsilon_j.
\]
For $i\ne j$, set
$\mb h_{ij}=\varepsilon_i\mb e_{u_i}-\varepsilon_j\mb e_{u_j}$.
By~\eqref{e-7},
\begin{gather*}
\mb h_{ij}^{\mathsf T}\left(R-\frac1\lambda I\right)\mb h_{ij}
=2p-2p=0,\\
(M_\alpha(H)\oplus0)\mb h_{ij}=0.
\end{gather*}
If $|\alpha|<1$, then~\eqref{e-3} would imply
\[
0\ge(1-|\alpha|)\bigl(d_H(u_i)+d_H(u_j)\bigr)>0.
\]
Hence $|\alpha|=1$ and $\mb g=\alpha\mb f$.
Since $e_i$ and $e_j$ meet only at their center outside $H$,
equality in~\eqref{e-13} implies
\[
s^2=\langle B_i\mb y_i,B_j\mb y_j\rangle
=(\mb g^{\mathsf T}\mb y_i)(\mb g^{\mathsf T}\mb y_j)
=\varepsilon_i\varepsilon_j s^2.
\]
Thus all $\varepsilon_i$ are equal, and the kernel relation above
reduces to $M_\alpha(H)(\mb e_{u_i}-\mb e_{u_j})=0$.
\end{proof}

\vspace{3mm}
\noindent\textbf{Claim 4.}
\emph{We have $\alpha=1$, $t=2$ and $\Delta(H)=1$.
Moreover, $u_1$ and $u_2$ are the endpoints of a
$K_2$ component of $H$.}

\begin{proof}[Proof of Claim 4.]
For any two distinct neighbors $u,v$ of the center, Claim~3 implies
\[
\begin{aligned}
0&=(\mb e_u-\mb e_v)^{\mathsf T}M_\alpha(H)(\mb e_u-\mb e_v)\\
 &=d_H(u)+d_H(v)-2\alpha A(H)_{uv}.
\end{aligned}
\]
Since $d_H(u),d_H(v)\ge1$, $A(H)_{uv}\in\{0,1\}$ and $|\alpha|\le1$,
\[
2\le d_H(u)+d_H(v)=2\alpha A(H)_{uv}\le2.
\]
Hence $\alpha=1$, $uv\in E(H)$ and $d_H(u)=d_H(v)=1$.
This holds for every pair of distinct vertices among $u_1,\ldots,u_t$,
so $t-1\le d_H(u_1)=1$. As $t\ge2$, we obtain $t=2$,
and $u_1u_2$ is a $K_2$ component of $H$.
Thus $G$ replaces this component by $K_3$, leaves the other
components unchanged, and may contain additional isolated vertices.

Since $q(K_2)=2$ and $q(K_3)=4$,
\[
\lambda=q(G)=\max\{4,q(H)\}>q(H),
\]
and therefore $\lambda=4$.
Choose $w\in V(H)$ with $d_H(w)=\Delta(H)$.
As $H_w^2$ contains $K_{1,\Delta(H)+2}$, Lemma~\ref{lem-2.3}
and subgraph monotonicity imply
\[
4=\lambda\ge q(K_{1,\Delta(H)+2})=\Delta(H)+3.
\]
Together with $e(H)\ge1$, this proves $\Delta(H)=1$.
\end{proof}

The resulting graph is precisely the configuration in
part~\textup{(ii)}.

\smallskip
\needspace{4\baselineskip}
\noindent\textbf{Sufficiency.}
Part~\textup{(i)} follows because isolated vertices do not
change $\rho_\alpha$.
For part~\textup{(ii)}, $H$ consists of $K_1$ and $K_2$ components,
with at least one $K_2$. By Lemma~\ref{lem-2.3},
\[
q(H_v^2)=
\begin{cases}
4,&d_H(v)=1,\\
3,&d_H(v)=0.
\end{cases}
\]
The graph $G$ consists of one $K_3$ and components isomorphic
to $K_1$ or $K_2$. Hence
\[
q(G)=4=\max_{v\in V(H)}q(H_v^2).
\]
\end{proof}

\section{The adjacency analogue}\label{sec-5}

Let $\rho(F)=\lambda_{\max}(A(F))$ be the adjacency
spectral radius of a graph $F$.
The comparison in Conjectures~\ref{conj-1.1} and~\ref{conj-1.2}
fails for the adjacency spectral radius, as the following example shows.

\begin{example}
Let $G=K_3$, and let $H$ be an induced copy of $K_2$ in $G$.
Then $t=e(G)-e(H)=2$, and $H_v^t\cong K_{1,3}$ for each
$v\in V(H)$.
Since $\rho(K_3)=2$ and $\rho(K_{1,3})=\sqrt3$, we obtain
\[
\rho(G)=2>\sqrt3=\max_{v\in V(H)}\rho(H_v^t).
\]
\end{example}

Theorem~\ref{thm-5.4} establishes the comparison when each component
of $G-V(H)$ has a unique attachment vertex in $H$ and induces
a bipartite graph together with that vertex.
We first prove the required resolvent estimate.

\begin{lemma}\label{lem-5.3}
Let $F$ be a connected bipartite graph with $s\ge1$ edges,
and let $r\in V(F)$ have degree $d$.
Order the vertices of $F$ with $r$ first, and write
\[
A(F)=
\begin{pmatrix}
0 & \mb b^{\mathsf T}\\
\mb b & C
\end{pmatrix},
\]
where $C=A(F-r)$.
For every real number $z>\sqrt{s}$,
\[
\mb b^{\mathsf T}(zI-C)^{-1}\mb b
\le\frac{zd}{z^2-s+d}\le\frac{s}{z}.
\]
Equality holds throughout if and only if $F$ is a star centered at $r$.
\end{lemma}

\begin{proof}
Let $U$ and $V$ be the parts of a bipartition of $F$, with
$r\in U$. Order the vertices of $F-r$ by listing
$U\setminus\{r\}$ before $V$. Then
\[
C=
\begin{pmatrix}
0&N\\
N^{\mathsf T}&0
\end{pmatrix},
\]
where $N$ is the biadjacency matrix of $F-r$ in this ordering.
The restriction $\mb b_0$ of $\mb b$ to $V$ records the
neighbors of $r$, while all other entries of $\mb b$ are zero.
Since $\mb b_0^{\mathsf T}\mb b_0=d$ and
$\operatorname{tr}(N^{\mathsf T}N)=s-d$, we have
\[
0\preceq N^{\mathsf T}N\preceq(s-d)I.
\]
In particular,
$\rho(F-r)^2=\lambda_{\max}(N^{\mathsf T}N)\le s-d<z^2$,
so $zI-C$ is positive definite. By block elimination and the
matrix bound above, we obtain
\[
\mb b^{\mathsf T}(zI-C)^{-1}\mb b
=z\,\mb b_0^{\mathsf T}
   (z^2I-N^{\mathsf T}N)^{-1}\mb b_0
\le\frac{z\,\mb b_0^{\mathsf T}\mb b_0}{z^2-s+d}
=\frac{zd}{z^2-s+d}.
\]
Also,
\begin{equation*}
\frac{s}{z}-\frac{zd}{z^2-s+d}
=
\frac{(s-d)(z^2-s)}
{z(z^2-s+d)}
\ge0.
\end{equation*}

Since $z^2>s$, equality in the second inequality holds
if and only if $d=s$. As $F$ is connected, this is equivalent
to $F$ being a star centered at $r$. In this case, $C=0$ and
$\mb b^{\mathsf T}\mb b=s$, so
$\mb b^{\mathsf T}(zI-C)^{-1}\mb b=s/z$.
\end{proof}

\begin{theorem}\label{thm-5.4}
Let $G$ be connected, and let $H$ be a connected induced
subgraph of $G$ with $e(H)\ge1$.
Suppose that $t=e(G)-e(H)\ge1$.
For every component $C$ of $G-V(H)$, assume that there is a
unique vertex $v_C\in V(H)$ adjacent to a vertex of $C$,
and that $G[V(C)\cup\{v_C\}]$ is bipartite.
Then
\[
\rho(G)\le\max_{v\in V(H)}\rho(H_v^t).
\]
Equality holds if and only if $G=H_v^t$ relative to the
prescribed copy of $H$, where $v$ attains the maximum
on the right.
\end{theorem}

\begin{proof}
Set $\lambda=\max_{v\in V(H)}\rho(H_v^t)$.
Attaching the $t$ new leaves to an endpoint of an edge of $H$
produces a graph containing $K_{1,t+1}$. By the variational principle,
\[
\lambda\ge\rho(K_{1,t+1})=\sqrt{t+1}>\sqrt t.
\]
For every $v\in V(H)$, we have $\lambda I-A(H_v^t)\succeq0$.
By Lemma~\ref{lem-2.2}, the Schur complement of the leaf block
$\lambda I_t$ satisfies
\begin{equation}\label{e-42}
T_v:=\lambda I-A(H)-\frac{t}{\lambda}\mb e_v\mb e_v^{\mathsf T}
\succeq0.
\end{equation}

Let $C_1,\ldots,C_a$ be the components of $G-V(H)$,
and let $v_i$ be the unique attachment vertex of $C_i$ in $H$.
For each $i$, set $F_i=G[V(C_i)\cup\{v_i\}]$ and $s_i=e(F_i)$,
and let $\mb b_i$ be the vector indexed by $V(C_i)$ whose entry
at $w$ is $1$ if $v_iw\in E(G)$ and $0$ otherwise.
Each $F_i$ is connected and bipartite, and
$E(G)\setminus E(H)=\bigcup_{i=1}^a E(F_i)$.
Hence $s_i\ge1$ and $\sum_{i=1}^a s_i=t$.
Figure~\ref{fig-3} illustrates this decomposition.

\begin{figure}[htb]
\centering
\begin{tikzpicture}[
  x=1.2cm,y=1cm,font=\small,
  vertex/.style={circle,fill=black,inner sep=0pt,minimum size=4pt},
  edge/.style={line width=0.65pt},
  extra/.style={edge,densely dashed},
  component/.style={draw=black!45,densely dotted,line width=0.5pt},
  line cap=round,line join=round]
  \filldraw[fill=black!5,draw=black!30,line width=0.4pt]
    (0.1,0) ellipse [x radius=1.25,y radius=1];
  \draw[component] (-2.6,0.45)
    ellipse [x radius=0.95,y radius=0.9];
  \draw[component] (2.5,0.45)
    ellipse [x radius=0.5,y radius=0.6];
  \coordinate (vone) at (-0.65,0.45);
  \coordinate (vtwo) at (0.7,0.45);
  \coordinate (h) at (0.05,-0.55);
  \coordinate (xone) at (-2.1,1.05);
  \coordinate (xtwo) at (-3.1,0.45);
  \coordinate (xthree) at (-2.1,-0.15);
  \coordinate (y) at (2.5,0.45);
  \draw[edge] (vone)--(vtwo)--(h)--(vone);
  \draw[extra] (vone)--(xone)--(xtwo)--(xthree)--(vone);
  \draw[extra] (vtwo)--(y);
  \foreach \u in {vone,vtwo,h,xone,xtwo,xthree,y}
    \node[vertex] at (\u) {};
  \node[above=3pt] at (vone) {$v_1$};
  \node[above=3pt] at (vtwo) {$v_2$};
  \node at (0.05,-0.05) {$H$};
  \node at (-2.6,1.1) {$C_1$};
  \node at (2.5,0.9) {$C_2$};
\end{tikzpicture}
\caption{The components $C_i$ and their attachment vertices $v_i$.}
\label{fig-3}
\end{figure}
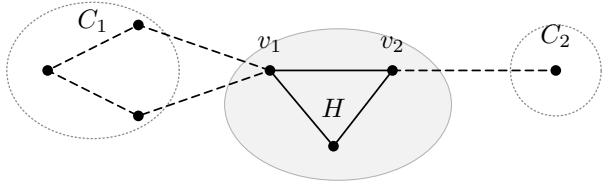

Since $C_i$ is bipartite,
$\rho(C_i)^2\le e(C_i)\le s_i\le t<\lambda^2$.
Thus $\lambda I-A(C_i)$ is positive definite, and by
Lemma~\ref{lem-5.3},
\[
\beta_i:=\mb b_i^{\mathsf T}(\lambda I-A(C_i))^{-1}\mb b_i
\le\frac{s_i}{\lambda}.
\]
The Schur complement of $\lambda I-A(G)$ on $V(H)$ is
\begin{equation}\label{e-43}
\begin{aligned}
S&:=\lambda I-A(H)-\sum_{i=1}^a\beta_i\mb e_{v_i}\mb e_{v_i}^{\mathsf T}\\
 &=\sum_{i=1}^a\frac{s_i}{t}T_{v_i}
   +\sum_{i=1}^a\left(\frac{s_i}{\lambda}-\beta_i\right)
       \mb e_{v_i}\mb e_{v_i}^{\mathsf T}\succeq0,
\end{aligned}
\end{equation}
where we used~\eqref{e-42} and the bound on $\beta_i$.
By Lemma~\ref{lem-2.2}, we obtain $\lambda I-A(G)\succeq0$,
and hence $\rho(G)\le\lambda$.

Suppose that $\rho(G)=\lambda$.
Let $\mb x>0$ satisfy $A(G)\mb x=\lambda\mb x$, and let
$\mb y=\mb x|_{V(H)}>0$.
Eliminating the coordinates outside $H$, we obtain $S\mb y=0$.
Since every summand in~\eqref{e-43} is positive semidefinite
and $s_i/t>0$, for each $1\le i\le a$ we obtain
$\left(\frac{s_i}{\lambda}-\beta_i\right)y_{v_i}^2=0$ and
$T_{v_i}\mb y=0$.
As $y_{v_i}>0$, we have $\beta_i=s_i/\lambda$.
By the equality condition in Lemma~\ref{lem-5.3}, each $F_i$
is a star centered at $v_i$.
Since $C_i=F_i-v_i$ is both independent and connected,
$|V(C_i)|=1$. Thus every vertex outside $H$ is pendant.
For any $i,j$, subtracting $T_{v_i}\mb y=0$ and
$T_{v_j}\mb y=0$, we obtain
\[
y_{v_i}\mb e_{v_i}=y_{v_j}\mb e_{v_j}.
\]
Since $\mb y>0$, all attachment vertices coincide.
Thus $G=H_v^t$ for some $v\in V(H)$, and
$\rho(H_v^t)=\rho(G)=\lambda$.
Conversely, equality is immediate for $H_v^t$ whenever
$\rho(H_v^t)=\max_{u\in V(H)}\rho(H_u^t)$.
\end{proof}

\medskip

\noindent\textbf{Acknowledgments.} The authors developed the arguments and wrote and verified all proofs. During an early exploratory stage, they used language-model-based tools to explore possible proof strategies. Suggestions from these tools served only as informal inspiration.

\medskip
\noindent\textbf{Data availability.} No datasets were generated or analyzed in this study.

\smallskip
\noindent\textbf{Conflict of interest.} The authors report no relevant financial or non-financial interests.

\end{document}